\documentclass[12pt]{article}
\usepackage{latexsym,amsmath,amssymb,color}
\usepackage{setspace}
\usepackage{amsthm}
\usepackage{eucal}
\usepackage{hyperref}
\usepackage[linesnumbered,lined,boxed,commentsnumbered,vlined]{algorithm2e}

\newcommand{\Integer}{\mathbb Z}

\newtheorem{Theorem}{Theorem}

\newtheorem{Corollary}{Corollary}
\theoremstyle{definition} 
\theoremstyle{definition} 
\theoremstyle{definition} 
\theoremstyle{definition} 
\theoremstyle{definition} 
\newenvironment{Proof}{\begin{trivlist} \item[] {\bf Proof.}}{\hfill $\Box$\end{trivlist}}

\newcommand{\sE}{\mathcal{E}}
\newcommand{\sF}{\mathcal{F}}
\newcommand{\sM}{\mathcal{M}}

\newcommand{\sP}{\mathcal{P}}
\newcommand{\sS}{\mathcal{S}}

\newcommand{\Ri}{{\sS}_1}
\newcommand{\Rii}{{\sS}_2}

\newcommand{\cmF}{\CMcal{F}}

\renewcommand{\geq}{\geqslant}
\renewcommand{\leq}{\leqslant}

\title{A Perfect One-Factorisation of $K_{56}$}

\author{David A.~Pike\thanks{Department of Mathematics and Statistics, Memorial University of Newfoundland, St.~John's, NL, Canada.  {\sl dapike@mun.ca}}}

\begin{document}

\date{\today}

\maketitle

\begin{abstract}
In 1963, Anton Kotzig conjectured that for each $n \geq 2$ the complete graph $K_{2n}$ has a perfect one-factorisation
({\em i.e.}, a decomposition into perfect matchings such that each pair of perfect matchings of the decomposition induces a Hamilton cycle).
We affirmatively settle the smallest unresolved case for this conjecture.
\end{abstract}

\vspace*{\baselineskip}

\vspace*{\baselineskip}

\section{Introduction}

A {\em $k$-factor} of a graph $G$ is a $k$-regular spanning subgraph of $G$ and
a {\em $k$-factorisation} of a $\Delta$-regular graph $G$ is a partition of $E(G)$ into $\frac{\Delta}{k}$ $k$-factors.
Thus a 1-factor (which is frequently also called a {\em perfect matching}) consists of a set of $\frac{|V(G)|}{2}$ pairwise non-adjacent edges,
and a 1-factorisation $\cmF = \{ \sF_1, \sF_2, \ldots, \sF_{\Delta} \}$ is a set of $\Delta$ disjoint 1-factors whose union is $E(G)$.

One-factorisations of the complete graph $K_{2n}$ have been known to exist since at least the 1890s when
F\'{e}lix Walecki was attributed with elegant means of partitioning the edges of complete graphs into Hamilton cycles~\cite{Lucas}.
A review of Walecki's decomposition techniques, presented with modern terminology and notation, is provided by Brian Alspach~\cite{Alspach2008}.

Note that although any Hamilton cycle of a graph of even order easily partitions into two 1-factors,
the converse is not necessarily the case.
In general it can only be said that the union of two disjoint 1-factors of a graph
results in a 2-factor which may itself consist of one or more cycles.

In the event that a 1-factorisation $\cmF = \{ \sF_1, \sF_2, \ldots, \sF_{\Delta} \}$ of a $\Delta$-regular graph has the particularly strong property
that $\sF_i \cup \sF_j$ is the edge set of Hamilton cycle for each $1 \leq i < j \leq \Delta$,
then the 1-factorisation is said to be a {\em perfect 1-factorisation}.
Anton Kotzig asked in 1963 whether, for each integer $n \geq 2$, the complete graph $K_{2n}$ has a perfect 1-factorisation~\cite{Kotzig1964}.
While this question has not yet been fully settled, it is known to have an affirmative answer whenever $n$ or $2n-1$ is prime
as well as when $2n$ is one of the following values that are mentioned in~\cite{AndersenHCD}:
16, 28, 36, 40, 50, 126, 170, 244, 344, 530, 730, 1332, 1370, 1850, 2198, 2810,
3126, 4490, 6860, 6890, 11450, 11882, 12168, 15626, 16808, 22202, 24390, 24650, 26570, 29792,
29930, 32042, 38810, 44522, 50654, 51530, 52442, 63002, 72362, 76730, 78126, 79508, 103824, 148878, 161052, 205380,
226982, 300764, 357912, 371294, 493040, 571788, 1092728, 1225044.
By April 2007, Ian Wanless had reported the following additional perfect one-factorisations at~\cite{WanlessWWW}:
1295030,
2248092,
2476100, 2685620,
3307950,
3442952, 4657464, 5735340, 6436344, 1030302, 2048384, 4330748, 6967872, 7880600, 9393932, 11089568, 11697084, 13651920, 15813252, 18191448, 19902512, 22665188.


Eric Seah published a survey article about perfect 1-factorisations and their properties in 1991~\cite{Seah1991},
and they are also mentioned by Walter Wallis in chapters of two books printed in 1992 and 1997
(see Section~8 of~\cite{WallisBB1992} and Chapter~16 of~\cite{Wallis1997}).
A more recent survey regarding perfect 1-factorisations, written by Alex Rosa, is forthcoming
and should be consulted for further details about their history and theoretical advances~\cite{Rosa2018}.
In the eleven years that have passed since the orders listed in the previous paragraph were published in~\cite{AndersenHCD,WanlessWWW},
only one new value has been confirmed,
namely 52, which was established by Adam Wolfe ten years ago~\cite{Wolfe2009}.
Using techniques described by Wolfe, in this paper we settle the smallest unresolved case of Kotzig's conjecture
by finding a perfect 1-factorisation of $K_{56}$.

\section{Main Results}

\begin{Theorem}
A perfect 1-factorisation of $K_{56}$ exists.
\end{Theorem}

\begin{Proof}
In his search for a perfect 1-factorisation of $K_{52}$, Wolfe showed how a pair of
starters in $\Integer_{2m-1}$ could be merged to produce an even starter for $\Integer_{4m-2}$.
In turn an even starter for $\Integer_{4m-2}$ yields a 1-factorisation for $K_{4m}$, which
can then be tested to assess whether it is perfect.

A {\em starter} in $\Integer_{2m-1}$ consists of a set $\sS$ of $m-1$ disjoint unordered pairs $\{x_i,y_i\} \subset \{0,1,\ldots,2m-2\}$
such that for each $d \in \{1,2,\ldots,2m-2\}$ one of the $m-1$ pairs $\{x_i,y_i\}$ satisfies
the criterion that either $x_i - y_i \equiv d$ (mod $2m-1$) or $y_i - x_i \equiv d$ (mod $2m-1$).
Note that any single element of $\{0,1,\ldots,2m-2\}$ is permitted to be absent from the $m-1$ pairs
(whereas it is often conventional to require 0 to be the missing element).
Our interest is in $m=14$.
In Table~\ref{Table-StarterPair} we present two starters for $\Integer_{27}$.
The missing element for starter $\Ri$ is 9, and for starter $\Rii$ element 20 is absent.

\begin{table}[b]
\begin{center}
\begin{tabular}{|c|p{61mm}|}
\hline
\multicolumn{2}{|c|}{$\Ri$} \\ \hline
High & \parbox{61mm}{
  $\{0,  1\}$,
  $\{7, 11\}$,
  $\{12, 17\}$,
  $\{20, 26\}$,
  $\{16, 25\}$,
  $\{8, 18\}$,
  $\{10, 22\}$} \\ \hline
Low & \parbox{61mm}{
$\{2,  4\}$,
$\{3,  6\}$,
$\{14, 21\}$,
$\{15, 23\}$,
$\{13, 24\}$,
$\{5, 19\}$} \\ \hline
\end{tabular}
\hspace*{0.6mm}
\begin{tabular}{|p{61mm}|c|}
\hline
\multicolumn{2}{|c|}{$\Rii$} \\ \hline
\parbox{61mm}{
   $\{1,  2\}$,
   $\{6, 10\}$,
   $\{16, 21\}$,
   $\{12, 18\}$,
   $\{7, 25\}$,
   $\{5, 15\}$,
   $\{8, 23\}$}
& Low \\ \hline
\parbox{61mm}{
$\{24, 26\}$,
$\{19, 22\}$,
$\{4, 11\}$,
$\{9, 17\}$,
$\{3, 14\}$,
$\{0, 13\}$}
& High \\ \hline
\end{tabular}
\end{center}
\caption{Two starters for $\Integer_{27}$}
\label{Table-StarterPair}
\end{table}

An {\em even starter} for $\Integer_{2t-2}$ consists of a set $\sE$ of $t-2$ disjoint unordered pairs $\{x_i,y_i\} \subset \{0,1,\ldots,2t-3\}$
such that
for each $d \in \{1,2,\ldots,2t-3\} \setminus \{t-1\}$ one of the $t-2$ pairs $\{x_i,y_i\}$ satisfies
the criterion that either $x_i - y_i \equiv d$ (mod $2t-2$) or $y_i - x_i \equiv d$ (mod $2t-2$).


As part of Wolfe's merging construction, the unordered pairs of elements that comprise each starter for $\Integer_{2m-1}$
are individually assigned a ``high'' or ``low'' designation
in such a way that the pair with difference $d \in \{1,2,\ldots,m-1\}$ from starter $\Ri$ is high (resp.\ low)
if and only if the pair with difference $d$ from starter $\Rii$ is low (resp.\ high).
Pseudocode to describe the merging operation that produces an even starter for $\Integer_{4m-2}$ is given in Algorithm~\ref{Pseudocode}.

\SetAlCapSty{}
\SetAlCapSkip{1em}
\IncMargin{1.1em}
\begin{algorithm}[!h]
\SetKwData{Left}{left}
\SetKwData{This}{this}
\SetKwData{Up}{up}
\SetKwFunction{Union}{Union}
\SetKwFunction{FindCompress}{FindCompress}
\SetKwInOut{Input}{Input}
\SetKwInOut{Output}{Output}

\Input{Two starters $\Ri$ and $\Rii$ for $\Integer_{2m-1}$, with corresponding high/low designations for each pair in $\Ri \cup \Rii$}
\Output{An even starter for $\Integer_{4m-2}$}
\BlankLine


$\sP \leftarrow \emptyset$

$a$ $\leftarrow$ the missing element of $\Ri$

$a \leftarrow a + 2m-1$

\For{$i \leftarrow 1$ \KwTo $2m-2$}
{
$x$ $\leftarrow$ $a$ (mod $2m-1$)

$y$ $\leftarrow$ the sole element of $\{0,1,\ldots,2m-2\}$ such that $\{x,y\} \in (\Ri \cup \Rii) \setminus \sP$

$d \leftarrow \min \big\{ x - y$ (mod $2m-1$), $y - x$ (mod $2m-1$)$\big\}$

$Y \leftarrow \{y, y+2m-1\}$

$\hat{y} \leftarrow$ the sole element of $Y \cap \big\{ a-d$ (mod $4m-2$), $a+d$ (mod $4m-2$)$\big\}$

  \uIf{$\{x,y\}$ \upshape{is low}}{
  $b \leftarrow \hat{y}$
  }
  \Else{
  $b \leftarrow$ the sole element of $Y \setminus \{\hat{y}\}$
  }

Output the pair $\{a,b\}$

$\sP \leftarrow \sP \cup \big\{ \{ x,y \} \big\}$

$a \leftarrow$ the sole element of $Y \setminus \{b\}$
}

\caption{Merging two starters to yield an even starter}
\label{Pseudocode}
\end{algorithm}
\DecMargin{1.1em}

When Algorithm~\ref{Pseudocode} is applied to the starters $\Ri$ and $\Rii$ and the high/low designations shown
in Table~\ref{Table-StarterPair}, the following even starter for $\Integer_{54}$
is obtained:

\begin{center}
$\{36,17\}$,
$\{44,12\}$,
$\{39,45\}$,
$\{18,35\}$,
$\{8,50\}$,
$\{23,15\}$,
$\{42,32\}$,
$\{5,46\}$,
$\{19,49\}$, \\
$\{22,37\}$,
$\{10,6\}$,
$\{33,30\}$,
$\{3,41\}$,
$\{14,21\}$,
$\{48,43\}$,
$\{16,52\}$,
$\{25,34\}$,
$\{7,38\}$, \\
$\{11,31\}$,
$\{4,2\}$,
$\{29,28\}$,
$\{1,27\}$,
$\{0,40\}$,
$\{13,24\}$,
$\{51,26\}$,
$\{53,20\}$
\end{center}

To construct a 1-factorisation of $K_{2t}$ from an even starter $\sE$ for $\Integer_{2t-2}$,
first observe that two elements of $\{0,1,\ldots,2t-3\}$, say $a$ and $b$, are not present among the $t-2$ pairs of $\sE$.
Let $\sF_0 = \sE \cup \big\{ \{a,\infty_1\}, \{b,\infty_2\} \big\}$.
For each $i \in \{1,2,\ldots,2t-3\}$ let $\sF_i = \sigma^i (\sF_0)$,
whereby the permutation $\sigma = (0, 1, 2, \ldots, 2t-3)(\infty_1)(\infty_2)$ is applied $i$ times to each element of each pair of $\sF_0$.
Also let $\sM = \big\{ \{x,x+t-1\} : x \in \{0,1,\ldots,t-2\} \big\}  \cup \big\{ \{\infty_1,\infty_2\}\big\}$.
The 1-factors $\sM$ and $\sF_0, \sF_1, \ldots, \sF_{2t-3}$ constitute a 1-factorisation of the complete graph of order $2t$
with vertex set $\{0,1,\ldots,2t-3\} \cup \{\infty_1,\infty_2\}$.
The 1-factorisation $\cmF = \{ \sM, \sF_0, \sF_1, \ldots, \sF_{2t-3} \}$ is not assured to be perfect,
but it is straightforward to test any given 1-factorisation to determine whether each pair of its 1-factors induces a Hamilton cycle.

For the even starter that we have presented for $\Integer_{54}$,
it is elements 9 and 47 that are absent from the pairs of $\sE$.
The resulting 1-factorisation of $K_{56}$, when assessed, was found to be perfect.
\end{Proof}

To briefly comment on the computational effort that was exerted in finding this perfect 1-factorisation,
we made use of a cluster of
IBM/Lenovo NeXtScale nx360 M4 nodes, each configured with
two Intel Xeon CPU E5-2650 v2  2.60GHz processors.
The search for a perfect 1-factorisation was distributed among 1024 tasks running in parallel,
one of them operating as a director and the rest as workers.
%
The worker that found the pair of starters shown in Table~\ref{Table-StarterPair} reported that it had
built and compared 7730443 pairs of starters (along with $2^{13}$ high/low combinations for each pair)
over the course of 33 days and 6 hours.  The 1022 other workers continued with their work until they
were terminated, at which time they had each been running for 43 days and 9 hours.
Hence our discovery of a perfect 1-factorisation of $K_{56}$ entailed the comparison of an estimated 10.3 billion pairs of starters.
To compare with Wolfe's discovery of a perfect 1-factorisation of $K_{52}$, he reported that
his search involved 7.494 billion pairs of starters~\cite{Wolfe2009}.
Note, however, that in our search for a perfect one-factorsation of $K_{56}$ a number of unsuccessful previous attempts,
some of them also lasting for a month, preceded the successful instance that is described above;
accounting details for these earlier attempts are not available.

Having found a perfect 1-factorisation for $K_{56}$, some consequential results follow.
For instance, Philip Laufer has proved that if $K_{2n}$ has a perfect 1-factorisation then so does the
complete bipartite graph $K_{2n-1,2n-1}$~\cite{Laufer1980}.
Under certain conditions the converse of Laufer's result is also known to hold~\cite{Sarwate1983,WanlessIhrig2005}.

\begin{Corollary}
A perfect 1-factorisation of $K_{55,55}$ exists.
\end{Corollary}

The existence of a perfect 1-factorisation of
$K_{m,m}$
implies the existence of a Latin square of order $m$ with no proper subsquares~\cite{Wanless1999}.
However, we note that the existence of such a Latin square of order 55
has already been established by Katherine Heinrich~\cite{Heinrich1980}.

Applications of perfect 1-factorisations also extend to coding theory.
Lihao Xu {\em et al.}\ have shown that a perfect 1-factorisation of $K_{2n}$ is equivalent to
a B-code $B_{2n-1}$, which is a type of Maximum Distance Separable (MDS) array code of distance 3~\cite{XBBW1999}.

\begin{Corollary}
A B-code $B_{55}$ exists.
\end{Corollary}

Applications of perfect 1-factorisations and their corresponding B-codes to the
design of RAID schemes for distributed storage are described in~\cite{DeoNanda2005,HuangBruck2016,NandaDeo2006}.
Some other codes and related structures arising from perfect 1-factorisations are discussed in~\cite{Wang2015}.

We conclude by observing that with the discovery of a perfect 1-factorisation of $K_{56}$, the smallest complete graph for which it is
unknown whether a perfect 1-factorisation exists becomes $K_{64}$.
Other unsettled orders up to 100 are 66, 70, 76, 78, 92, 96 and 100.

\section{Acknowledgements}
The author acknowledges research grant support from NSERC
as well as computational support from
The Centre for Health Informatics and Analytics of the Faculty of Medicine at Memorial University of Newfoundland.



\begin{thebibliography}{99}

\bibitem{Alspach2008}B.\ Alspach.
The wonderful Walecki construction.
{\em Bull.\ Inst.\ Combin.\ Appl.} 52 (2008) 7--20.

\bibitem{AndersenHCD}L.D.\ Andersen.
Factorizations of Graphs,
in {\em The Handbook of Combinatorial Designs}, second edition
(eds.\ C.J.\ Colbourn and J.H.\ Dinitz).
Chapman and Hall / CRC, Boca Raton (2007) pp.\ 740--755.

\bibitem{DeoNanda2005}N.\ Deo and S.\ Nanda.
One-factors and Hamiltonian paths in modeling data and parity placement in disk arrays.
36th Southeastern International Conference on Combinatorics, Graph Theory, and Computing.
{\em Congr.\ Numer.} 176 (2005) 191--199.



\bibitem{Heinrich1980}K.\ Heinrich.
Latin squares with no proper subsquares.
{\em J.\ Combin.\ Theory Ser.\ A} 29 (1980) 346--353.


\bibitem{HuangBruck2016}W.\ Huang and J.\ Bruck.
Secure RAID schemes for distributed storage,
IEEE International Symposium on Information Theory
(2016) pp.\ 1401--1405.

\bibitem{Kotzig1964}A.\ Kotzig.
         Hamilton graphs and Hamilton circuits, in
         {\em Theory of Graphs and its Applications\/}, Proceedings of the Symposium held in Smolenice in June 1963.
         Academic Press, New York (1964) pp.\ 63--82 and 162.

\bibitem{Laufer1980}P.J.\ Laufer.
On strongly Hamiltonian complete bipartite graphs.
{\em Ars Combin.} 9 (1980) 43--46.



\bibitem{Lucas}\'{E}.\ Lucas,
Les jeux de demoiselles,
{\em R\'{e}cr\'{e}at.\ Math.\/} Vol.\ II,
        Gauthier-Villars et Fils, Paris (1892), 161--197.

\bibitem{NandaDeo2006}S.\ Nanda and N.\ Deo.
Methods for placing data and parity to tolerate two disk failures in disk arrays using complete bipartite graphs.
37th Southeastern International Conference on Combinatorics, Graph Theory, and Computing.
{\em Congr.\ Numer.} 179 (2006) 167--179.

\bibitem{Rosa2018}A.\ Rosa.
Perfect 1-factorizations.
{\em Math.\ Slovaca}, in press.

\bibitem{Sarwate1983}D.G.\ Sarwate.
A note on strongly regular Hamiltonian equivalence of $K_{2n}$ and $K_{2n-1,2n-1}$ and its generalization.
{\em  J. Indian Acad. Math.} 5 (1983) 65--67.

\bibitem{Seah1991}E.\ Seah.
Perfect one-factorizations of the complete graph -- a survey.
{\em Bull.\ Inst.\ Combin.\ Appl.} 1 (1991) 59--70.

\bibitem{WallisBB1992}W.D.\ Wallis.
One-factorizations of complete graphs,
in {\em Contemporary Design Theory: A Collection of Surveys}
(eds.\ J.H.\  Dinitz and D.R.\ Stinson). Wiley, New York (1992) pp.\ 593--631.

\bibitem{Wallis1997}W.D.\ Wallis.
{\em One-Factorizations}.
Kluwer, Dordrecht (1997).

\bibitem{Wang2015}Y.\ Wang.
Privacy-preserving data storage in cloud using array BP-XOR codes.
{\em  IEEE Trans.\ Cloud Comput.} 3 (2015) 425--435.


\bibitem{Wanless1999}I.M.\ Wanless.
Perfect factorisations of bipartite graphs and Latin squares without proper subrectangles.
{\em Electron.\ J.\ Combin.} 6 (1999) \#R9.


\bibitem{WanlessWWW}I.M.\ Wanless.
\url{http://users.monash.edu.au/~iwanless/data/P1F/newP1F.html}


\bibitem{WanlessIhrig2005}I.M.\ Wanless and E.C.\ Ihrig.
Symmetries that Latin squares inherit from 1-factorizations.
{\em J.\ Combin.\ Des.} 13 (2005) 157--172.


\bibitem{Wolfe2009}A.J.\ Wolfe.
A perfect one-factorization of $K_{52}$.
{\em J.\ Combin.\ Des.} 17 (2009) 190--196.


\bibitem{XBBW1999}L.\ Xu, V.\ Bohossian, J.\ Bruck and D.G.\ Wagner.
Low-density MDS codes and factors of complete graphs.
{\em IEEE Trans.\ Inform.\ Theory} 45 (1999) 1817--1826.

\end{thebibliography}
\end{document}